\documentclass[11pt]{article}
\usepackage{amsfonts}

\usepackage{graphicx}
\usepackage{latexsym}
\usepackage{amsmath}
\usepackage{layout}

\newtheorem{theorem}{Theorem}
\newtheorem{remark}{Remark}
\newtheorem{example}{Example}

\begin{document}

\title{Occupation densities for certain   processes related to fractional Brownian motion }
\author{ Khalifa Es-Sebaiy $^{1}\quad$ David Nualart $^{2}\quad$ Youssef
Ouknine $^{1}\quad$ Ciprian A. Tudor $^{3}\vspace*{0.1in}$ \\
%EndAName
$^{1}$ Department of Mathematics, Faculty of Sciences Semlalia,\\
Cadi Ayyad University 2390 Marrakesh, Morocco. \vspace*{0.1in}\\
$^{2}$ Department of Mathematics, University of Kansas,\\
405 Snow Hall, Lawrence, Kansas 66045-2142, USA,\vspace*{0.1in}\\
$^{3}$SAMOS/MATISSE, Centre d'Economie de La Sorbonne,\\
Universit\'e de Panth\'eon-Sorbonne Paris 1,\\
90, rue de Tolbiac, 75634 Paris Cedex 13, France.}
\maketitle

\begin{abstract}
In this paper we establish the existence of a square integrable occupation density for two classes of stochastic processes.
First we consider a Gaussian process with an absolutely continuous random drift, and secondly we handle the case of a
(Skorohod) integral with respect to the fractional Brownian motion with Hurst parameter $H>\frac 12$.
 The proof of these results  uses a general criterion for the existence of a square integrable local time, which is based on the techniques of Malliavin calculus.
\end{abstract}

\vskip0.5cm


 \vskip0.2cm


\renewcommand{\thefootnote}{\fnsymbol{footnote}}

\vskip1cm

\section{Introduction}

Local times for semimartingales have been widely studied. See for example the monograph \cite{RY} and the references therein.
On the other hand, local times of Gaussian processes have also been the object of a rich probabilistic literature; see for
example the recent paper \cite{MR} by Marcus and Rosen. A general criterion for the existence of a   local time for a wide
class of anticipating processes, which are not semimartingales or Gaussian processes,  was established by Imkeller and Nualart
in \cite{NuIm}. The proof of this result combines the techniques of
Malliavin calculus with the criterion given by Geman and Horowitz in \cite%
{GH}. This criterion was applied in \cite{NuIm} to the Brownian motion with
an anticipating drift, and to indefinite Skorohod integral processes.

The aim of this paper is to establish the existence of the occupations
densities for two classes of stochastic processes related to the fractional Brownian motion, using the approach introduced in %
\cite{NuIm}. First we consider a Gaussian   process $B=\{B_{t}, t\in [0,1]\}$ with  an absolutely continuous random drift
\begin{equation*}
X_{t}= B_{t} +\int_{0}^{t} u_{s}ds,
\end{equation*}
where $u$ is a stochastic process measurable with respect to the $\sigma$%
-field generated by $B$.  We assume
that the variance of the increment of the Gaussian process $B$ on an
interval $[s,t]$ behaves as $|t-s|^{2\rho}$, for some $\rho\in (0,1)$. This
includes, for instance, the bifractional Brownian motion with parameters $%
H,K\in (0,1)$.
Under reasonable regularity hypotheses imposed to
the process $u$ we prove the existence of a square integrable occupation
density with respect to the Lebesque measure for the process $X$.

Our second example is represented by the indefinite divergence (Skorohod)
integral $X=\{X_{t}, t\in [(0,1]\}$ with respect to the fractional Brownian
motion with Hurst parameter $H\in (\frac 12, 1)$, that is
\begin{equation*}
X_{t}= \int_{0}^{t} u_{s} \delta B^H_{s}.
\end{equation*}
We provide  integrability conditions on the integrand $u$ and its  iterated  derivatives in the sense of Malliavin calculus in
order to deduce the existence of a square integrable occupation densities for $X$.

We organized our paper as follows. Section 2 contains some preliminaries on
the Malliavin calculus with respect to Gaussian processes. In Section 3 we
prove the existence of the occupation densities for perturbed Gaussian
processes and in Section 4 we treat the case of indefinite divergence
integral processes with respect to the fractional Brownian motion.

\vskip0.5cm

\section{Preliminaries}

Let $\{B_{t} , t\in [0,1]\}$ be a centered Gaussian process with covariance
function
\begin{equation*}
R(t,s):=E( B_{t}B_{s}),
\end{equation*}
defined in a complete probability space $(\Omega, \mathcal{F},P)$. By ${%
\mathcal{H}}$ we denote the canonical Hilbert space associated to $B$ defined as the closure of the linear space generated by
the indicator functions $\{ \mathbf{1}_{[0,t]}, t\in [0,1]\} $ with respect to the inner product
\begin{equation*}
\langle \mathbf{1}_{[0,t]} , \mathbf{1}_{[0,s] } \rangle _{{\mathcal{H}}}
=R(t,s), \hskip0.5cm s,t\in [0,1].
\end{equation*}
The mapping $\mathbf{1}_{[0,t]} \to X_{t}$ can be extended to an isometry
between ${\mathcal{H}}$ and the first Gaussian chaos generated by $B$. We
denote by $B(\varphi) $ the image of an element $\varphi \in {\mathcal{H}}$
by this isometry.

We will first introduce some elements of the Malliavin calculus associated
with $B$. We refer to \cite{N06} for a detailed account of these notions.
For a smooth random variable $F=f\left( B(\varphi _{1}), \ldots ,
B(\varphi_{n})\right) $, with $\varphi_{i} \in {\mathcal{H}}$ and $f\in
C_{b}^{\infty}(R^{n})$ ($f$ and all its partial derivatives are bounded) the
derivative of $F$ with respect to $B$ is defined by
\begin{equation*}
D F =\sum_{j=1}^{n}\frac{\partial f}{\partial x_{j}}(B(\varphi_{1}),\dots,B(%
\varphi_{n}))\varphi_{j}.
\end{equation*}
For any integer $k\ge 1$ and any real number $p\ge 1$ we denote by $\mathbb{D%
}^{k,p}$ the Sobolev space defined as the the closure of the space of smooth
random variables with respect to the norm
\begin{equation*}
\Vert F\Vert_{k,p}^{p}=E(|F|^{p})+\sum_{j=1}^{k}\Vert
D^{j}F\Vert_{L^{p}(\Omega;{\mathcal{H}}^{\otimes j})}^{p}.
\end{equation*}
Similarly, for a given Hilbert space $V$ we can define Sobolev spaces of $V$%
-valued random variables $\mathbb{D}^{k,p}(V)$.

Consider the adjoint $\delta $ of $D $ in $L^2$. Its domain is the class of
elements $u\in L^{2}(\Omega;{\mathcal{H}})$ such that
\begin{equation*}
E(\langle D F,u\rangle_{{\mathcal{H}}})\leq C\Vert F\Vert_{2},
\end{equation*}
for any $F\in \mathbb{D} ^{1,2}$, and $\delta \left( u\right) $ is the
element of $L^{2}(\Omega)$ given by
\begin{equation*}
E(\delta (u)F)=E(\langle D F,u\rangle_{{\mathcal{H}}})
\end{equation*}
for any $F\in \mathbb{D}^{1,2}$. We will make use of the notation $\delta
(u)=\int_{0}^{1}u_{s}\delta B_{s}$. It is well-known that $\mathbb{D}^{1,2}({%
\mathcal{H}})$ is included in the domain of $\delta $. Note that $E(\delta (
u ) )=0$ and the variance of $\delta(u)$ is given by
\begin{equation}
E(\delta (u)^{2})=E(\Vert u\Vert_{{\mathcal{H}}}^{2})+E(\langle D u,(D
u)^{\ast}\rangle_{{\mathcal{H}}\otimes{\mathcal{H}}} ),  \label{squaremean}
\end{equation}
if $u\in \mathbb{D}^{1,2}({\mathcal{H}})$, where $(D u)^{\ast}$ is the
adjoint of $D u$ in the Hilbert space ${\mathcal{H}}\otimes{\mathcal{H}}$.
We have Meyer's inequality
\begin{equation}
E(|\delta (u)^{p}|)\leq C_{p}\left( E(\Vert u\Vert_{{\mathcal{H}}%
}^{p})+E(\Vert D u\Vert_{{\mathcal{H}}\otimes{\mathcal{H}}}^{p})\right),
\label{meyer}
\end{equation}
for any $p>1$. We will make use of the property
\begin{equation}
F\delta (u)=\delta (Fu)+\langle D F,u\rangle_{{\mathcal{H}}}.  \label{p1}
\end{equation}
if $F\in \mathbb{D}^{1,2}$ and $u\in \mathrm{Dom}(\delta )$ such that $Fu\in
\mathrm{Dom}(\delta )$. We also need the commutativity relationship between $%
D $ and $\delta $
\begin{equation}  \label{comm}
D \delta (u)= u + \int_{0}^{1} D u_{s}\delta B_{s},
\end{equation}
if $u\in \mathbb{D}^{1,2}({\mathcal{H}})$ and the process $\{D_{s}u, s\in
[0,1]\}$ belongs to the domain of $\delta $.

Throughout this paper we will assume that the centered Gaussian process $%
B=\{B_{t},t\in \lbrack 0,1]\}$ satisfies
\begin{equation}
C_{1}(t-s)^{2\rho }\leq E(|B_{t}-B_{s}|^{2})\leq C_{2}(t-s)^{2\rho },
\label{cond-cov}
\end{equation}%
for some $\rho\in (0,1)$ with $C_{1},C_{2}$ two positive constants not
depending on $t,s$. It will follow from the Kolmogorov criterium that $B$
admits a H\"{o}lder continuous version of order $\delta $ for any $\delta
<\rho $.

Throughout this paper we will denote by $C$ a generic constant that may be
different from line to line.

\begin{example}
The bifractional Brownian motion (see, for instance \cite{HV}), denoted by $B^{H,K}$, is defined as a
centered Gaussian process starting from zero with covariance
\begin{equation}  \label{covbiFBM}
R(t,s)= \frac{1}{2^{K}}\left( \left( t^{2H} + s^{2H}\right) ^{K} -\vert
t-s\vert ^{2HK}\right)
\end{equation}
where $H\in (0,1) $ and $K\in (0,1]$. When $K=1$, then we have a standard fractional Brownian motion denoted by $B^H$. It has
been proven in \cite{HV} that for all $s\leq t$,
\begin{equation}  \label{qh}
2^{-K}\vert t-s\vert ^{2HK} \leq E\left| B^{H,K}_{t}- B^{H,K}_{s}\right|
^{2}\leq 2^{1-K}\vert t-s\vert ^{2HK}
\end{equation}
so relation (\ref{cond-cov}) holds with $\rho=HK$. A stochastic analysis for this process can be found in \cite{KRT} and a
study of its occupation densities  has been done in  \cite{ET},  \cite{XT}.
\end{example}

  For a measurable function $x: [0,1]\to \mathbb{R}$ we define the
occupation measure
\begin{equation*}
\mu (x)(C) =\int_0^1 \mathbf{1}_{C}(x_{s})ds,
\end{equation*}
where $C $ is a Borel subset of $\mathbb{R}$ and we will say that $x$ has on
occupation density with respect to the Lebesque measure $\lambda $ if the
measure $\mu $ is absolutely continuous with respect to $\lambda$. The
occupation density of the function $x$ will be the derivative $\frac{d\mu
_{t}}{d\lambda }$. For a continuous process $\{X_{t}, t\in [0,1]\}$ we will
say that $X$ has an occupation density on $[0,1]$ if for almost all $\omega
\in \Omega$, $X(\omega)$ has an occupation density on $[0,1]$.

We will use the following criterium for the existence of occupation
densities (see \cite{NuIm}). Set $T=\{(s,t)\in \lbrack 0,1]^{2}:s<t\}$.

\begin{theorem}
\label{th1} Let $\{X_{t},t\in \lbrack 0,1]\}$ be a continuous stochastic
process such that $X_{t}\in \mathbb{D}^{2,2}$ for every $t\in \lbrack 0,1]$.
Suppose that there exists a sequence of random variables $\{F_{n},n\geq 1\}$
with $\bigcup_{n}\{F_{n}\not=0\}=\Omega $ a.s. and $F_{n}\in \mathbb{D}%
^{1,1} $ for every $n\geq 1$, two sequences $\alpha _{n}>0,\delta _{n}>0$,
a measurable bounded function $\gamma :[0,1]\rightarrow \mathbb{R}$, and a
constant $\ \theta >0$, such that:
\begin{description}
\item[a)] For every $n\geq 1$, $|t-s|\leq \delta _{n}$, and on $%
\{F_{n}\not=0\}$ we have
\begin{equation}
\langle \gamma D(X_{t}-X_{s}), \mathbf{1}_{[s,t]}\rangle _{{\mathcal{H}}%
}>\alpha _{n}|t-s|^{\theta },\qquad \mathrm{a.s.}.  \label{a}
\end{equation}
\item[b)] For every $n\geq 1$
\begin{equation}
\int_{T}E(\langle \gamma DF_{n},\mathbf{1}_{[s,t]}\rangle _{\mathcal{H}%
})|t-s|^{-\theta }dtds<\infty.  \label{b}
\end{equation}
\item[c)] For every $n\geq 1$
\begin{equation}
\int_{T}E\left( \left| F_{n}\left\langle \gamma ^{\otimes 2}DD(X_{t}-X_{s}),%
\mathbf{1}_{[s,t]}^{\otimes 2}\right\rangle _{\mathcal{H}^{\otimes 2}}\right| \right) |t-s|^{-2\theta }dsdt<\infty .
\label{c}
\end{equation}
\end{description}
Then the process $\{X_{t},t\in \lbrack 0,1]\}$ admits a square integrable
occupation density on $[0,1]$.
\end{theorem}

\begin{remark}
The original result has been stated in \cite{NuIm} with $\theta =1$ in the case of the standard Brownian motion. On the other
hand, by applying Proposition 2.3 and Theorem 2.1 in \cite{NuIm} it follows easily that  this criterium can be stated for  any
$\theta
>0$.
\end{remark}

\setcounter{equation}0

\vskip0.5cm

\section{Occupation density for Gaussian processes  with random drift}

We study in this part the existence of the occupation density for Gaussian
processes perturbed by a absolute continuous random drift. The main result
of this section is the following.

\begin{theorem}
Let $\{B_{t},t\in \lbrack 0,1]\}$ be a Gaussian process satisfying (\ref%
{cond-cov}). Consider the process $\{X_{t},t\in \lbrack 0,1]\}$ given by
\begin{equation*}
X_{t}=B_{t}+\int_{0}^{t}u_{s}ds,
\end{equation*}%
and suppose that the process $u$ satisfies the following conditions:
\begin{enumerate}
\item  $u \in \mathbb{D}^{2,2}(L^{2}([0,1]))$.
\item $E\left( \left( \int_{0}^{1}\left\| D^{2}u_{t}\right\| _{\mathcal{%
H\otimes H}}^{p}dt\right) ^{q/p}\right) <\infty $, for some $q>1$, $p>\frac{1%
}{1-\rho }$.
\end{enumerate}
Then, the process $X$ has a square integrable occupation density on the
interval $[0,1]$.
\end{theorem}

\vskip10pt\noindent \textit{Proof:}\hskip10pt We are going to apply Theorem %
\ref{th1}. Notice first that $X_{t}\in \mathbb{D}^{2,2}$ for all $t\in
\lbrack 0,1]$. For any $0\leq s<t\leq 1$, using (\ref{comm}) and (\ref%
{cond-cov}) we have%
\begin{eqnarray*}
\left\langle D(X_{t}-X_{s}),\mathbf{1}_{[s,t]}\right\rangle _{\mathcal{H}}
&=&\left\langle \mathbf{1}_{[s,t]},\mathbf{1}_{[s,t]}\right\rangle _{%
\mathcal{H}}+\left\langle \int_{s}^{t}Du_{r}dr,\mathbf{1}_{[s,t]}\right%
\rangle _{\mathcal{H}} \\
&\geq &C_{1}(t-s)^{2\rho }-\left| \left\langle \int_{s}^{t}Du_{r}dr,\mathbf{1%
}_{[s,t]}\right\rangle _{\mathcal{H}}\right|  \\
&\geq &C_{1}(t-s)^{2\rho }-\sqrt{C_{2}}(t-s)^{\rho }\int_{s}^{t}\left\|
Du_{r}\right\| _{\mathcal{H}}dr.
\end{eqnarray*}%
By H\"{o}lder's inequality, if $\frac{1}{p}+\frac{1}{q}=1$, we obtain%
\begin{equation*}
\int_{s}^{t}\left\| Du_{r}\right\| _{\mathcal{H}}dr\leq (t-s)^{\frac{1}{q}%
}\left( \int_{0}^{1}\left\| Du_{r}\right\| _{\mathcal{H}}^{p}dr\right) ^{%
\frac{1}{p}}.
\end{equation*}%
Fix a natural number $n\geq 2$, and choose a function $\varphi _{n}(x)$,
which is infinitely differentiable with compact support, such that $\varphi
_{n}(x)=1$ if $|x|\leq n-1$, and $\varphi _{n}(x)=0$, if $|x|\geq n$. Set $%
F_{n}=\varphi _{n}\left( \left( \int_{0}^{1}\left\| Du_{t}\right\| _{%
\mathcal{H}}^{p}dt\right) ^{\frac{1}{p}}\right) $. The random variable $F_{n}
$ belogs to $\mathbb{D}^{1,q}$. In fact, it suffices to write $F_{n}=\varphi
_{n}(G)$, where%
\begin{equation*}
G=\sup_{\substack{ h\in L^{q}([0,1];\mathcal{H)} \\ \left\| h\right\| \leq 1
}}\int_{0}^{1}\left\langle Du_{r},h_{r}\right\rangle _{\mathcal{H}}dr,
\end{equation*}%
which implies%
\begin{eqnarray*}
\left\| DF_{n}\right\| _{\mathcal{H}} &=&\left\| \varphi _{n}^{\prime
}(G)DG\right\| _{\mathcal{H}}\leq \left\| \varphi _{n}^{\prime }\right\|
_{\infty }\sup_{\substack{ h\in L^{q}([0,1];\mathcal{H)} \\ \left\|
h\right\| \leq 1}}\left\| \int_{0}^{1}\left\langle
D^{2}u_{r},h_{r}\right\rangle _{\mathcal{H}^{\otimes 2}}dr\right\| _{%
\mathcal{H}} \\
&\leq &\left\| \varphi _{n}^{\prime }\right\| _{\infty }\left(
\int_{0}^{1}\left\| D^{2}u_{r}\right\| _{\mathcal{H}^{\otimes
2}}^{p}dr\right) ^{\frac{1}{p}}\in L^{q}(\Omega ).
\end{eqnarray*}%
Then, on \ the set $\{F_{n}\neq 0\}$, $\left( \int_{0}^{1}\left\|
Du_{t}\right\| _{\mathcal{H}}^{p}dt\right) ^{\frac{1}{p}}\leq n$, and we get%
\begin{eqnarray*}
\left\langle D(X_{t}-X_{s}),\mathbf{1}_{[s,t]}\right\rangle _{\mathcal{H}}
&\geq &C_{1}(t-s)^{2\rho }-n\sqrt{C_{2}}(t-s)^{\rho +\frac{1}{q}} \\
&=&(t-s)^{2\rho }\left[ C_{1}-n\sqrt{C_{2}}(t-s)^{\frac{1}{q}-\rho }\right] ,
\end{eqnarray*}%
and \ property a) of Theorem \ref{th1} holds with a suitable choice of $%
\alpha _{n}$ and $\delta _{n}$ because $\frac{1}{q}-\rho >0$, and with $%
\theta =2\rho $ and $\gamma =1$.

Finally, conditions b) and c) can also be checked:%
\begin{equation*}
\int_{T}\frac{E\left( \left| \left\langle DF_{n},\mathbf{1}%
_{[s,t]}\right\rangle _{\mathcal{H}}\right| \right) }{|t-s|^{2\rho }}%
dsdt\leq \sqrt{C_{2}}\int_{T}\frac{E\left( \left\| DF_{n}\right\| _{\mathcal{%
H}}\right) }{|t-s|^{\rho }}dsdt<\infty ,
\end{equation*}%
and%
\begin{eqnarray*}
&&\int_{T}\frac{E\left( \left| F_{n}\left\langle D^{2}(X_{t}-X_{s}),\mathbf{1%
}_{[s,t]}^{\otimes 2}\right\rangle _{\mathcal{H}^{\otimes 2}}\right| \right)
}{|t-s|^{4\rho }}dsdt \\
&\leq &\left\| F_{n}\right\| _{\infty }C_{2}\int_{T}\frac{E\left( \left\|
D^{2}(X_{t}-X_{s})\right\| _{\mathcal{H}^{\otimes 2}}\right) }{|t-s|^{2\rho }%
}dsdt<\infty ,
\end{eqnarray*}%
because%
\begin{eqnarray*}
E\left( \left\| D^{2}(X_{t}-X_{s})\right\| _{\mathcal{H}^{\otimes 2}}\right)
&=&E\left( \left\| \int_{s}^{t}D^{2}u_{r}dr\right\| _{\mathcal{H}^{\otimes
2}}\right) \leq \int_{s}^{t}E\left( \left\| D^{2}u_{r}\right\| _{\mathcal{H}%
^{\otimes 2}}\right) dr \\
&\leq &(t-s)^{\frac{1}{q}}E\left[ \left( \int_{0}^{1}\left\|
D^{2}u_{r}\right\| _{\mathcal{H}^{\otimes 2}}^{p}dr\right) ^{\frac{1}{p}}%
\right] ,
\end{eqnarray*}%
and $\frac{1}{q}-2\rho =1-\frac{1}{p}-2\rho >-1$, because \ $p>\frac{1}{%
2(1-\rho )}$. \hfill \vrule width.25cm height.25cm depth0cm\smallskip

\vskip0.5cm

\begin{remark}
These conditions are intrinsic and they do not depend on the structure of the Hilbert space $\mathcal{H}$. In the case of the
Brownian motion, this result is slightly weaker than Theorem 3.1 in \cite{NuIm}, because we require a little more
integrability.
\end{remark}

\vskip0.5cm

\setcounter{equation}0
\section{Occupation density for  Skorohod integrals with respect to the
fractional Brownian motion}

We study here the existence of occupation densities for indefinite
divergence integrals with respect to the fractional Brownian motion.
Consider a process of the form $X_{t}=\int_{0}^{t}u_{s}\delta B_{s}^{H}$, $%
t\in \lbrack 0,1]$, where $B$ is fractional Brownian motion with Hurst
parameter $H\in \left( \frac{1}{2},1\right) $, and $u$ is an element of $%
\mathbb{D}^{1,2}(L^{2}([0,1]))\subset \mathrm{Dom}\left( \delta \right) $.

We know that the covariance of the fractional Brownian motion can be written
as
\begin{equation}
E(B_{t}^{H}B_{s}^{H})=\int_{0}^{t}\int_{0}^{s}\phi (\alpha ,\beta )d\alpha
d\beta ,  \label{H1}
\end{equation}%
where $\phi (\alpha ,\beta )=H(2H-1)|\alpha -\beta |^{2H-2}$. For any $0\leq
s<t\leq 1$, and $\alpha \in \lbrack 0,1]$ we set
\begin{equation}
f_{s,t}(\alpha ):=\int_{s}^{t}\phi (\alpha ,\beta )d\alpha d\beta .
\label{f}
\end{equation}
We also know (see e.g. \cite{N06}) that the canonical Hilbert space associated to $B$ satisfies:
\begin{equation}
\label{inclu} L^{2}\left( [0,1]\right) \subset L^{\frac{1}{H}}\left( [0,1]\right)\subset \mathcal{H}.
\end{equation}

The following is the main result of this section.

\begin{theorem}
Consider the stochastic process $X_{t}=\int_{0}^{t}u_{s}\delta B_{s}^{H}$
where the integrand $u$ satisfy the following conditions for some $q>\frac{2H%
}{1-H}$ and $p>1$ such that $\frac 1p +2 < H(p+1)$:
\begin{description}
\item[I1)] {\ }$u\in \mathbb{D}^{3,2}(L^{2}([0,1]))$.
\item[I2)] {\ }$\int_{0}^{1} \int_0^1 [ E( | D_{t}u_s |^p) +E( \| \left| D_{t}Du_s\right|
\|_{\mathcal{H}} ^p) +E( \| \left| D_{t}DDu_s\right|  \|_{\mathcal{H}\otimes \mathcal{H}} ^p) ] dsdt  <\infty $.
\item[I3)] $\int_{0}^{1}E\left( |u_{t}|^{-\frac{p}{p-1}(q+1)}\right)
dt<\infty $.
\end{description}
Then the process $\{X_{t},t\in \lbrack 0,1]\}$ admits a square integrable
occupation density on $[0,1]$.
\end{theorem}

\vskip10pt\noindent \textit{Proof:}\hskip10pt We will use the criteria given
in \cite{NuIm} and recalled in Theorem \ref{th1}. Condition I1) implies that
$X_{t}\in \mathbb{D}^{2,2}$ for all $t\in \lbrack 0,1]$. On the other hand,
from Theorem 7.8 in \cite{KRT} (or also by a slightly modification of
Theorem 5 in \cite{AN}) we obtain the continuity of the paths of the process
$X$. Note that from Lemma 2.2 in \cite{NuIm} corroborated with hypothesis
I3). we obtain the existence of a function $\gamma :[0,1]\rightarrow \{-1,1\}
$ such that $\gamma _{t}u_{t}=|u_{t}|$ for almost all $t$ and $\omega $.

We are going to show conditions a), b) and c) of Theorem \ref{th1}.

\vskip 0.3cm
\noindent  \textit{Proof of condition a): } Fix $0\leq s<t\leq 1$. From (\ref%
{comm}) we obtain
\begin{equation*}
D(X_{t}-X_{s})=u\mathbf{1}_{[s,t]}+\int_{s}^{t}Du_{r}\delta B_{r}^{H},
\end{equation*}%
and we can write
\begin{equation}
\langle \gamma (X_{t}-X_{s}), \mathbf{\ 1}_{[s,t]}\rangle _{\mathcal{H}%
}=\langle |u|\mathbf{1}_{[s,t]},\mathbf{1}_{[s,t]}\rangle _{\mathcal{H}%
}+\langle \gamma \int_{s}^{t}Du_{r}\delta B_{r}^{H},\mathbf{1}%
_{[s,t]}\rangle _{\mathcal{H}}.  \label{z1}
\end{equation}%
We first study the term
\begin{equation*}
\langle |u|\mathbf{1}_{[s,t]},\mathbf{1}_{[s,t]}\rangle _{\mathcal{H}%
}=\int_{s}^{t}\int_{s}^{t}|u_{\alpha }|\phi (\alpha ,\beta )d\alpha d\beta
=\int_{s}^{t}|u_{\alpha }|f_{s,t}(\alpha )d\alpha .
\end{equation*}%
For any $q>1$ we have
\begin{eqnarray*}
E(|B_{t}^{H}-B_{s}^{H}|^{2}) &=&\int_{s}^{t}f_{s,t}(\alpha )d\alpha  \\
&=&\int_{s}^{t}\left( |u_{\alpha }|f_{s,t}(\alpha )\right) ^{\frac{q}{q+1}%
}\left( |u_{\alpha }|f_{s,t}(\alpha )\right) ^{-\frac{q}{q+1}}f_{s,t}(\alpha
)d\alpha ,
\end{eqnarray*}%
and using H\"{o}lder's inequality with orders $\frac{q+1}{q}$ and $q+1$, we
obtain
\begin{equation*}
E(|B_{t}^{H}-B_{s}^{H}|^{2})\leq \left( \int_{s}^{t}|u_{\alpha
}|f_{s,t}(\alpha )d\alpha \right) ^{\frac{q}{q+1}}\left(
\int_{s}^{t}|u_{\alpha }|^{-q}f_{s,t}(\alpha )d\alpha \right) ^{\frac{1}{q+1}%
}.
\end{equation*}%
Hence, using that%
\begin{equation*}
f_{s,t}(\alpha )\leq f_{0,1}(\alpha )=H(2H-1)\int_{0}^{1}|\alpha -\beta |^{2H-2}d\beta =H\left( \alpha ^{2H-1}+(1-\alpha
)^{2H-1}\right) \leq H,
\end{equation*}%
we get
\begin{equation}
\int_{s}^{t}|u_{\alpha }|f_{s,t}(\alpha )d\alpha \geq C |t-s|^{\frac{2H(q+1)}{q%
}}Z_{q}^{-\frac{1}{q}},  \label{z2}
\end{equation}%
where $Z_{q}=$ $\int_{0}^{1}|u_{\alpha }|^{-q}d\alpha $.

On the other hand, for the second summand in the right-hand side of (\ref{z1}%
) we can write, using H\"{o}lder's inequality.
\begin{eqnarray}
\left| \left\langle \gamma  \int_{s}^{t}Du_{r}\delta B_{r}^{H},\mathbf{1}%
_{[s,t]}\right\rangle _{\mathcal{H}}\right|  &\leq &\int_{0}^{1}\left| \int_{s}^{t}D_{\alpha }u_{r}\delta B_{r}^{H}\right|
f_{s,t}(\alpha )d\alpha
\notag \\
&\leq &\left( \int_{0}^{1}f_{s,t}(\alpha )^{\frac{p}{p-1}}d\alpha \right) ^{%
\frac{p-1}{p}}  \notag \\
&&\times \left( \int_{0}^{1}\left| \int_{s}^{t}D_{\alpha }u_{r}\delta
B_{r}\right| ^{p}d\alpha \right) ^{\frac{1}{p}}.  \label{bb2}
\end{eqnarray}%
We can write%
\begin{eqnarray}
\left( \int_{0}^{1}f_{s,t}(\alpha )^{\frac{p}{p-1}}d\alpha \right) ^{\frac{%
p-1}{p}} &=&c_{H}\left\| \int_{s}^{t}|\cdot -\beta |^{2H-2}d\beta \right\|
_{L^{\frac{p}{p-1}}([0,1])}  \notag \\
&\leq &c_{H}\left\| \mathbf{1}_{[s,t]}\ast |\cdot |^{2H-2}\mathbf{1}%
_{[-1,1]}\right\| _{L^{\frac{p}{p-1}}(\mathbb{R})},  \label{bb3}
\end{eqnarray}%
where $c_{H}=H(2H-1)$. Young's inequality with exponents $a$ and $b$ in $%
(1,\infty )$ such that $\frac{1}{a}+\frac{1}{b}=2-\frac{1}{p}$ yields
\begin{equation}
\left\| \mathbf{1}_{[s,t]}\ast |\cdot |^{2H-2}\mathbf{1}_{[-1,1]}\right\|
_{L^{\frac{p}{p-1}}(\mathbb{R})}\leq \left\| \mathbf{1}_{[s,t]}\right\|
_{L^{a}(\mathbb{R})}\left\| |\cdot |^{2H-2}\mathbf{1}_{[-1,1]}\right\|
_{L^{b}(\mathbb{R})}.  \label{bb4}
\end{equation}%
Choosing $b<\frac{1}{2-2H}$ and letting $\eta =\frac{1}{a}<2H-\frac{1}{p}$ we obtain from (\ref{bb2}), (\ref{bb3}), and
(\ref{bb4})
\begin{equation*}
\left| \left\langle \gamma \int_{s}^{t}Du_{r}\delta B_{r}^{H},\mathbf{1}%
_{[s,t]}\right\rangle _{\mathcal{H}}\right| \leq C|t-s|^{\eta }\left( \int_{0}^{1}\left| \int_{s}^{t}D_{\alpha }u_{r}\delta
B_{r}\right| ^{p}d\alpha \right) ^{\frac{1}{p}}.
\end{equation*}%
Now we will apply Garsia-Rodemich-Ramsey's lemma (see \cite{GRR}) with $\Phi
(x)=x^{p}$, $p(x)=x^{\frac{m+2}{p}}$ and to the continuous function $%
u_{s}=\int_{0}^{s}D_{\alpha }u_{r}\delta B_{r}$ (use again Theorem 5 in \cite%
{AN}), and we get
\begin{equation}
\left| \left\langle \int_{s}^{t}Du_{r}\delta B_{r},\gamma \mathbf{1}%
_{[s,t]}\right\rangle _{\mathcal{H}}\right| \leq C|t-s|^{\eta +\frac{m}{p}}%
 Y_{m,p}^{\frac{1}{p}},  \label{bb5}
\end{equation}%
where%
\begin{equation*}
Y_{m,p}=\int_{0}^{1}\int_{0}^{1}\int_{0}^{1}\frac{\left|
\int_{x}^{y}D_{\alpha }u_{r}\delta B_{r}\right| ^{p}}{|x-y|^{m+2}}%
dxdyd\alpha .
\end{equation*}%
Substituting (\ref{z2}) and (\ref{bb5}) into (\ref{z1}) yields
\begin{eqnarray*}
\langle \gamma D(X_{t}-X_{s}),1_{[s,t]}\rangle _{\mathcal{H}} &\geq &|t-s|^{%
\frac{2H(q+1)}{q}}Z_{q}^{-\frac{1}{q}}-C|t-s|^{\eta +\frac{m}{p}}Y_{m,p}^{%
\frac{1}{p}} \\
&=&|t-s|^{\frac{2H(q+1)}{q}}\left( Z_{q}^{-\frac{1}{q}}-C|t-s|^{\delta
}Y_{m,p}^{\frac{1}{p}}\right) ,
\end{eqnarray*}%
where $\delta =\eta +\frac{m}{p}-2H-\frac{2H}{q}$. With a right choice of $%
\eta $ the exponent $\delta $ is positive, provided that $m-\frac{1}{p}-%
\frac{2H}{q}>0$, because $\eta <2H-\frac{1}{p}$. \ Taking into account that
 $\frac{2H}{q}<1-H$, it suffices that%
\begin{equation}
m>\frac{1}{p}+1-H.  \label{c1}
\end{equation}

We construct now the sequence $\left\{ F_{n},n\geq 1\right\} $. \ Fix a
natural number $n\geq 2$, and choose a function $\varphi _{n}(x)$, which is
infinitely differentiable with compact support, such that $\varphi _{n}(x)=1$
if $|x|\leq n-1$, and $\varphi _{n}(x)=0$, if $|x|\geq n$. Set $%
F_{n}=\varphi _{n}\left( G\ \right) $, where $G=Z_{q}+Y_{m,p}$. Then clearly  the sequences $\alpha _{n}$ and $\delta _{n}$
required in Theorem \ref{th1} can be constructed on the set $\left\{ F_{n}\not=0\right\} $, with $\theta =2H+\frac{2H}{q}$.

It only remains to show that the random variables $F_{n}$ are in the space $%
\mathbb{D}^{1,1}$. \ For this we have to show that the random variables \ $%
\left\| DZ_{q}\right\| _{\mathcal{H}}$ and $\left\| DY_{m,p}\right\| _{%
\mathcal{H}}$ are integrable on the set $\{G\leq n\}$. First notice that, as
in the proof of \ Proposition  4.1   of \cite{NuIm}, we can show that  $%
E\left( \left\| DZ_{q}\right\| _{\mathcal{H}}\right) <\infty $.  This
follows from the integrability conditions I3) and%
\begin{equation}
\int_{0}^{1}E\left( \Vert Du_{t}\Vert _{\mathcal{H}}^{p}\right) dt<\infty ,
\label{c3}
\end{equation}%
which holds because of \ I2), the continuous embedding of of $L^{\frac{1}{H}%
}([0,1])$ into $\mathcal{H}$ (see \cite{MMV}), and the fact that $pH\geq 1$.
On the other hand, we can write
\begin{equation*}
DY_{m,p}=p\int_{0}^{1}\int_{0}^{1}\int_{0}^{1}\left| \xi _{x,y,\alpha
}\right| ^{p-1}\mathrm{sign}(\xi _{x,y,\alpha })D\xi _{x,y,\alpha
}|x-y|^{-m-2}dxdyd\alpha ,
\end{equation*}%
where $\xi _{x,y,\alpha }=\int_{y}^{x}D_{\alpha }u_{r}\delta B_{r}$. Thus
\begin{eqnarray*}
\Vert DY_{m,p}\Vert _{\mathcal{H}} &\leq
&p\int_{0}^{1}\int_{0}^{1}\int_{0}^{1}\left| \xi _{x,y,\alpha }\right|
^{p-1}\Vert D\xi _{x,y,\alpha }\Vert _{\mathcal{H}}|x-y|^{-m-2}dxdyd\alpha
\\
&\leq &p(Y_{m,p})^{\frac{p-1}{p}}\left(
\int_{0}^{1}\int_{0}^{1}\int_{0}^{1}\Vert D\xi _{x,y,\alpha }\Vert _{%
\mathcal{H}}^{p}|x-y|^{-m-2}dxdyd\alpha \right) ^{1/p}.
\end{eqnarray*}%
Now, to show that $1_{(G\leq n)}\Vert DY_{m,p}\Vert _{\mathcal{H}}$ belongs
to $L^{1}(\Omega )$, it suffices to show that the random variable
\begin{equation*}
Y=\int_{0}^{1}\int_{0}^{1}\int_{0}^{1}\Vert D\xi _{x,y,\alpha }\Vert _{%
\mathcal{H}}^{p}|x-y|^{-m-2}dxdyd\alpha
\end{equation*}%
has a finite expectation. \  Since, for any $0\leq y<x\leq 1$
\begin{equation*}
D\xi _{x,y,\alpha }=\mathbf{1}_{[y,x]}D_{\alpha }u+\int_{y}^{x}DD_{\alpha
}u_{s}\delta B_{s}^{H},
\end{equation*}%
we have
\begin{eqnarray*}
Y &\leq &C\left( \int_{0}^{1}\int_{0}^{1}\int_{0}^{1}\Vert \mathbf{1}%
_{[y,x]}D_{\alpha }u\Vert _{\mathcal{H}}^{p}|x-y|^{-m-2}dxdyd\alpha \right.
\\
&&+\left. \int_{0}^{1}\int_{0}^{1}\int_{0}^{1}\Vert \int_{y}^{x}DD_{\alpha
}u_{s}\delta B_{s}^{H}\Vert _{\mathcal{H}}^{p}|x-y|^{-m-2}dxdyd\alpha
\right)  \\
:= &&C(Y_{1}+Y_{2}).
\end{eqnarray*}%

From the  continuous embedding of $L^{\frac{1}{H}}([0,1])$ into $\mathcal{H}$%
, we obtain%
\begin{eqnarray*}
Y_{1} &\leq &C\int_{0}^{1}\int_{0}^{1}\int_{0}^{1}\Vert \mathbf{1}%
_{[y,x]}D_{\alpha }u\Vert _{L^{1/H}([0,1])}^{p}|x-y|^{-m-2}dxdyd\alpha  \\
&\leq &C|x-y|^{pH-1}\int_{0}^{1}\int_{0}^{1}\int_{0}^{1}\int_{y}^{x}\left|
D_{\alpha }u_{r}\right| ^{p}|x-y|^{-m-2}drdxdyd\alpha .
\end{eqnarray*}%
Hence, $E(Y_{1})<\infty $, by Fubini's theorem, Proposition 3.1 in \cite{NuIm} and condition I2), provided%
\begin{equation}
m<pH-1.  \label{c1a}
\end{equation}
On the other hand, using the estimate (\ref{meyer}), and again  the
continuous embedding of $L^{\frac{1}{H}}([0,1])$ into $\mathcal{H}$ yields%
\begin{eqnarray*}
E\left( \Vert \int_{y}^{x}DD_{\alpha }u_{s}\delta B_{s}^{H}\Vert _{\mathcal{H%
}}^{p}\right)  &\leq &C\ E\left( \left\| D_{\alpha }Du_{\cdot}\mathbf{1}%
_{[y,x]}(\cdot)\right\| _{\mathcal{H}^{\otimes 2}}^{p}+\left\| D_{\alpha }DDu_{\cdot}\mathbf{1}%
_{[y,x]}(\cdot)\right\| _{\mathcal{H}^{\otimes 3}}^{p}\right)  \\
&\leq &C\ E\Big( \left\| \left|D_{\alpha }Du_{\cdot }\right| \mathbf{1}_{[y,x]}(\cdot)\right\|
_{L^{1/H}([0,1];\mathcal{H})}^{p} \\
&&+\left\| \left|  D_{\alpha }DDu_{\cdot} \right| \mathbf{1}%
_{[y,x]}(\cdot )\right\| _{L^{1/H}([0,1];\mathcal{H}^{\otimes 2})}^{p}\Big)  \\
&\leq &C|x-y|^{pH-1}\Bigg ( \int_{y}^{x}E\left( \left\| \left|  D_{\alpha
}Du_{r}\right| \right\| _{\mathcal{H}}^{p}\right) dr  \\
&& +\int_{y}^{x}E\left( \left\| \left| D_{\alpha }DDu_{r}\right| \right\| _{\mathcal{H}^{\otimes 2}}^{p}\right) dr\Bigg ) .
\end{eqnarray*}%
As before we obtain $E(Y_{2})<\infty $  by Fubini's theorem and condition
I2), provided (\ref{c1a}) holds.     \ Notice that \ condition \ $\frac{1}{p}%
+2<H(p+1)$ implies that we can choose an $m$ such that \ (\ref{c1}) and (\ref%
{c1a}) hold.

\vskip0.3cm  \noindent \textit{Proof of condition b): } Define $A_{n}=\{G\leq n\}$.
Then, condition b) in Theorem \ref{th1} follows from
\begin{eqnarray*}
\int_{T}E(\langle \gamma DF_{n},\mathbf{1}_{[s,t]}\rangle _{\mathcal{H}%
})|t-s|^{-\theta }dtds &\leq &C\int_{T}E(\mathbf{1}_{A_{n}}\left| \langle \gamma DG, \mathbf{1}_{[s,t]}\rangle
_{\mathcal{H}}\right| )|t-s|^{-\theta
}dtds \\
&\leq &CE\left( \mathbf{1}_{A_{n}}\Vert DG\Vert _{\mathcal{H}}\right)
\int_{T}|t-s|^{H-\theta }dsdt<\infty ,
\end{eqnarray*}%
since $E\left( \mathbf{1}_{A_{n}}\Vert DG\Vert _{\mathcal{H}}\right) <\infty
$ and $\theta -H=H+\frac{2H}{q}<1$.

\vskip0.3cm \noindent  \textit{Proof of condition c):} We have
\begin{equation*}
D_{\alpha }D_{\beta }(X_{t}-X_{s})=\mathbf{1}_{[s,t]}(\beta )D_{\alpha
}u_{\beta }+\mathbf{1}_{[s,t]}(\alpha )D_{\beta }u_{\alpha
}+\int_{s}^{t}D_{\alpha }D_{\beta }u_{r}\delta B_{r}^{H}.
\end{equation*}%
Hence$\ $
\begin{eqnarray*}
\left\langle \gamma ^{\otimes 2} DD(X_{t}-X_{s}),\mathbf{1}_{[s,t]}^{\otimes
2}\right\rangle _{\mathcal{H}^{\otimes 2}} &=&\left\langle \gamma ^{\otimes 2}\mathbf{1}%
_{[s,t]}(\beta )D_{\alpha }u_{\beta },\mathbf{1}%
_{[s,t]}^{\otimes 2}\right\rangle _{\mathcal{H}^{\otimes 2}}+\left\langle \gamma ^{\otimes 2}
\mathbf{1}_{[s,t]}(\alpha )D_{\beta }u_{\alpha },\mathbf{1%
}_{[s,t]}^{\otimes 2}\right\rangle _{\mathcal{H}^{\otimes 2}} \\
&&+\left\langle \gamma ^{\otimes 2}\int_{s}^{t}D_{\alpha }D_{\beta }u_{r}\delta
B_{r}^{H},\mathbf{1}_{[s,t]}^{\otimes 2}\right\rangle _{%
\mathcal{H}^{\otimes 2}} \\
:= &&J_{1}(s,t)+J_{2}(s,t)+J_{3}(s,t).
\end{eqnarray*}%
For $i=1,2,3$, we set%
\begin{equation*}
A_{i}=E\left( F_{n}\int_{T}|t-s|^{-2\theta }\left| J_{i}(s,t)\right|
dsdt\right) .
\end{equation*}%
Let us compute first%
\begin{eqnarray*}
A_{1} &\leq &C\int _{T}\int_{T}|t-s|^{2H-2\theta }E\left( \Vert |D_{\alpha }u_{\beta
}|\mathbf{1}_{[s,t]}(\beta )\Vert _{\mathcal{H}^{\otimes 2}}\right) dsdt \\
&=&C\ \int_{T}\int_{T}|t-s|^{2H-2\theta }\left( \int_{s}^{t}\int_{s}^{t}\varphi (\beta ,y)d\beta dy\right)^{\frac{1}{2}} dsdt,
\end{eqnarray*}%
where%
\begin{equation*}
\varphi (\beta ,y)=\int_{0}^{1}\int_{0}^{1}E\left( |D_{\alpha }u_{\beta
}||D_{x}u_{y}|\right) \phi (\alpha ,x)\phi (\beta ,y)d\alpha dx.
\end{equation*}%
By Fubini's theorem $A_{1}<\infty $, because $2H-2\theta >-2$, which is
equivalent to $q>H$, and%
\begin{equation*}
\int_{0}^{1}\int_{0}^{1}\varphi (\beta ,y)d\beta dy \leq E\left( \left\| \left| Du \right| \right\| _{%
\mathcal{H}^{\otimes 2}}^{2}\right)
\end{equation*}%
and this is finite because of the inclusion of $L^{2}([0,1])$ in $\mathcal{H}$ (\ref{inclu}). In the same way we can show that
$A_{2}<\infty $. Finally,
\begin{eqnarray*}
A_{3} &=&E\left( F_{n}\int_{T}\int_{T}|t-s|^{-2\theta }\left| \left\langle \gamma ^{\otimes 2}
\int_{s}^{t}D_{\alpha }D_{\beta }u_{r}\delta B_{r}^{H},%
\mathbf{1}_{[s,t]}^{\otimes 2}\right\rangle _{\mathcal{H}^{\otimes
2}}\right| dsdt\right)  \\
&\leq &C\int_{T}\int_{T}|t-s|^{2H-2\theta }E\left( \left\| \int_{s}^{t}DDu_{r}\delta B_{r}^{H}\right\| _{\mathcal{H}^{\otimes
2}}\right) dsdt,
\end{eqnarray*}%
and we conclude as before by using for example the bound (\ref{meyer}) for the norm of the Skorohod integral and the condition
I2). \hfill \vrule width.25cm height.25cm depth0cm\smallskip

\vskip0.5cm

\begin{remark}
If $p= \frac {1+\sqrt{17}}2$, then $\frac 1p +2 < H(p+1)$ for all $H>\frac 12
$.
\end{remark}

\end{document}